\documentstyle[twoside,12pt]{article}

\newcommand{\chapter}{\section}

\begin{document}

\newtheorem{Thm}{Theorem}
\newtheorem{Ax}{Axiom}
\newtheorem{Prop}{Proposition}
\newtheorem{Cor}[Prop]{Corollary}
\newtheorem{Main}{}
\renewcommand{\theMain}{}
\newtheorem{Lem}[Prop]{Lemma}
\newtheorem{Fact}{Fact}
\renewcommand{\theFact}{}

\newtheorem{Def}{Definition}
\newtheorem{rmk}{Remark}
\newenvironment{Rmk}{\begin{rmk}\em}{\end{rmk}}
\newtheorem{exm}{Example}
\newenvironment{Exm}{\begin{exm}\em}{\end{exm}}

\newcommand{\qed}{\par {\EM QED} }
\newtheorem{prf}{Proof}
\renewcommand{\theprf}{}
\newenvironment{Prf}{\begin{prf}\em}{\qed\end{prf}}
\newtheorem{prff}{}
\renewcommand{\theprff}{}
\newenvironment{Prff}{\begin{prff}\em}{\qed\end{prff}}





\newcommand{\YES}[1]{#1}
\newcommand{\NOT}[1]{}

\newcommand{\cA}{{\cal A}}
\newcommand{\cB}{{\cal B}}
\newcommand{\cC}{{\cal C}}
\newcommand{\cD}{{\cal D}}
\newcommand{\cE}{{\cal E}}
\newcommand{\cF}{{\cal F}}
\newcommand{\cG}{{\cal G}}
\newcommand{\cH}{{\cal H}}
\newcommand{\cI}{{\cal I}}
\newcommand{\cJ}{{\cal J}}
\newcommand{\cK}{{\cal K}}
\newcommand{\cL}{{\cal L}}
\newcommand{\cM}{{\cal M}}
\newcommand{\cN}{{\cal N}}
\newcommand{\cO}{{\cal O}}
\newcommand{\cP}{{\cal P}}
\newcommand{\cQ}{{\cal Q}}
\newcommand{\cR}{{\cal R}}
\newcommand{\cS}{{\cal S}}
\newcommand{\cT}{{\cal T}}
\newcommand{\cU}{{\cal U}}
\newcommand{\cV}{{\cal V}}
\newcommand{\cW}{{\cal W}}
\newcommand{\cX}{{\cal X}}
\newcommand{\cY}{{\cal Y}}
\newcommand{\cZ}{{\cal Z}}

\newcommand{\bbb}[1]{{\mbox{\bf #1}}}

\newcommand{\bN}{\bbb{N}}
\newcommand{\bZ}{\bbb{Z}}
\newcommand{\bR}{\bbb{R}}
\newcommand{\bC}{\bbb{C}}
\newcommand{\bQ}{\bbb{Q}}
\newcommand{\bT}{\bbb{T}}

\newcommand{\noind}[1]{{\setlength{\parindent}{0cm} #1}}
\newcommand{\parsk}{\par\medskip}

\newcommand{\varend}{\end{document}}

\newcommand{\LP}{\left(}  \newcommand{\RP}{\right)}
\newcommand{\LQ}{\left[}  \newcommand{\RQ}{\right]}
\newcommand{\LB}{\left\{} \newcommand{\RB}{\right\}}
\newcommand{\LA}{\left\langle} \newcommand{\RA}{\right\rangle}
\newcommand{\LS}{\left|}  \newcommand{\RS}{\right|}
\newcommand{\LN}{\left\|} \newcommand{\RN}{\right\|}
\newcommand{\bLP}{\bigl(}  \newcommand{\bRP}{\bigr)}

\newcommand{\ts}{{\otimes}} \newcommand{\TS}{{\bigotimes}}

\newcommand{\es}{\emptyset}
\newcommand{\sm}{\setminus} \newcommand{\st}{\subset}
\newcommand{\eps}{\varepsilon}

\newcommand{\beware}{~$\Psi\!\Psi\!\Psi$~}

\newcommand{\EM}{\em\bf}

\newcommand{\BE}{\begin{equation}}  \newcommand{\EE}{\end{equation}}
\newcommand{\BER}{$$\begin{array}}  \newcommand{\EER}{\end{array}$$}
\newcommand{\head}[1]{
             \par\bigskip\noind{{\large\bf#1}\nopagebreak\par}\medskip}
\newcommand{\chead}[1]{\begin{center}
              \par\bigskip{\large\bf#1}\nopagebreak\par\medskip
              \end{center}}

\newcommand{\fillitem}{\L{\embox{}}\vspace{-.6cm}}

\newcommand{\dfrac}[2]{\displaystyle{\frac{#1}{#2}}}
\newcommand{\toprove}[1]{{\buildrel\mbox{?}\over#1}}
\newcommand{\DEF}{{\buildrel\mbox{\scriptsize\,\,def\,\,}\over=}}
\newcommand{\BAR}[1]{\overline{#1}}
\newcommand{\I}{\infty}
\newcommand{\al}{\alpha} \newcommand{\be}{\beta}
\newcommand{\OM}{\Omega} \newcommand{\om}{\omega}
\newcommand{\la}{\lambda}
\newcommand{\mod}{\mbox{ mod }}
\newcommand{\half}{\frac{1}{2}}
\newcommand{\NI}{{n\to\infty}}
\newcommand{\ang}{{<\!\!\!)}}
\newcommand{\arc}[1]{\breve{#1}}

\newenvironment{txeqn}[1]{\begin{equation}\label{#1}\begin{array}{l}}{
                         \end{array}\end{equation}}



\newcommand{\Middle}{{\mbox{Middle}\,}}
\newcommand{\End}{{\mbox{End}\,}}
\newcommand{\GL}{{\mbox{GL}\,}}
\newcommand{\diam}{{\mbox{diam}\,}}
\newcommand{\tr}{{\mbox{tr}\,}}

\title{How to Axiomatize School Geometry}
\author{Eliahu Levy\\
Department of Mathematics\\
Technion -- Israel Institute of Technology,
Haifa 32000, Israel\\
email: eliahu@techunix.technion.ac.il}


\date{}


\maketitle



This is an attempt to present axioms for Euclidean Geometry, aiming at the
following goals:
\begin{itemize}
 \item To work with ``geometrical'' notions. Thus we would not merely
identify points in the plane with pairs of real numbers, which means that a particular
coordinate system is given special status.
 \item To be appropriate to the way geometry is done in higher mathematics
(including physics and engineering). This means that the algebraic nature
of much of geometry need not be hidden.
 \item To respond to the desire that one would confidently accept
empirically/intuitively
that the axioms are valid in our physical everyday world (or rather in the
usual idealization that ``geometry'' is). This seems to disfavor
taking the Theorem of Pythagoras as an axiom.
 \item To have accessible the usual rigor of ``pure'' mathematics, and to
make the axioms satisfying by the standards of the latter. In particular,
not to take as an axiom something that can be naturally proved. Note that some
``topological'' notions, necessary for the rigor of the presented axioms, can
be readily kept silent with an ``unsophisticated'' audience
(such as school).
\end{itemize}
\medskip

The style in the sequel is intended for those accustomed to mathematical
writings, in order to make the mathematical contents clear. Of course, in case
an approach in this spirit can be practiced in school the style of
presentation must be quite different.


\section{The Axioms: Plane Geometry}

Primitive notions:

\begin{itemize}
\item A set $P$ (The {\EM Plane}), whose elements are called {\EM points}.
For the sophisticated -- the Plane is assumed a Hausdorff
topological space. We shall let $x,y,z,\;a,b,c$ etc.\ vary over
points.
\item A relation among 3 points, indeed a commutative ``algebraic'' operation:
$z$ is the middle between $x$ and $y$ (to be written $z=\Middle(x,y)$).
\item An equivalence relation on $P\times P$: two pairs of points have the same
distance.
\end{itemize}
\bigskip

We shall have altogether three axioms.
\medskip

By $\bR$ we denote, as usual, the set of {\EM real numbers}
(for the unsophisticated - just the set of {\EM numbers}, representable, say,
as possibly unending decimal fractions).

{\Ax \label{Ax:Co}
{\EM (Axiom of Coordinates)}. There is a bijection
(coordinates system) between $\bR^2$ and the Plane (which is a homeomorphism)
such that any mapping given in the coordinates
by $x\mapsto a+x$ or $x\mapsto a-x$ (here $a\in\bR^2\,,x\in\bR^2$) ``preserves
the geometry'': it is an ``isomorphism'' with respect to the middle operation
and maps any pair of points to a pair with the same distance.
}
\bigskip

Note that we did not take the coordinates as a primitive notion: there may
be many such bijections, the axiom saying that there is at least one.
\medskip

The empirical/intuitive evidence for this axiom is plain: one encounters such
coordinates daily (with the ``rough'' everyday correlate of the idealized
``set of points''). Maybe it is more intuitive to postulate the stronger
requirements that the reflections $(x_1,x_2)\mapsto(\alpha\pm x_1,x_2)$ and
$(x_1,x_2)\mapsto(x_1,\alpha\pm x_2)$\quad($\alpha\in\bR$) preserve the
geometry, the existence of enough reflections being implicit even in Euclid.
\medskip

As a simple consequence we may prove

{\Thm \label{T} In any coordinates system satisfying Axiom \ref{Ax:Co}
the middle operation $\Middle(a,b)$ corresponds to the ``algebraic middle''
$\dfrac{a+b}2\quad a,b\in\bR^2$.}

\begin{Prf}
The map $x\mapsto(a+b)-x$ maps $a\to b$, $b\to a$, and preserves the
geometry. Hence it fixes $\Middle(a,b)$. But its only fixed point is
$\dfrac{a+b}2$.
\end{Prf}
\medskip

Define, for integer $n\ge2$, an {\EM $n$-ruler} as a
sequence of points $(a_i)_{0\le i\le n}$ such that for any $0<i<n$\quad
$a_i=\Middle(a_{i-1},a_{i+1})$
\medskip
\NOT{
In fact, we need only $2$-rulers and the following Fact might be stated
only for $n=2$.
Introducing $n$-rulers for general $n$ serves only for didactic purposes.
\medskip}

Theorem \ref{T} implies that for any coordinates system as in Axiom
\ref{Ax:Co}, a ruler is just an ``algebraic ruler'', i.e.\ a sequence
with constant difference. We may deduce:
\medskip

{\Fact For any integers $n\ge k>l\ge 0$ and points $a,b$ there is
a unique $n$-ruler $(c_i)_{0\le i\le n}$ with $a=c_k$ and $b=c_l$.
}
\medskip

If we define, with respect to some coordinates system, an
``algebraic {\EM straight line}'' as usual (as a set $L\subset\bR^2$ of the form
$L=\{a+\lambda c\,|\,\lambda\in\bR\}$ where $a,c\in\bR^2,\,c\ne0$) then the
straight line joining $a,b\in P$ contains the ``rational line'' joining
$a$ and $b$, i.e.\ the set of all points obtained by constructing $n$-rulers
according to the Fact, and is its closure. Thus the notion of straight line
is independent of the coordinates (note that we needed the topology
here, and that just $2$-rulers would have sufficed).
\medskip

A quadrangle $(a,b,c,d)$ is an ``algebraic parallelogram'' with respect to
some coordinates system if $a-b=d-c$. But this is equivalent to $(a,c)$
and $(b,d)$ having the same ``algebraic middle'', i.e.\ to
$\Middle(a,c)=\Middle(b,d)$. Thus the notion of parallelogram is again
``geometrical'' -- independent of the coordinates system.
This allows us to define the {\EM vectors} geometrically as ``differences of
pairs of points'', that is, say, as equivalence classes of pairs of points
by the equivalence relation defined by parallelograms. (Thus a point minus
a point is a vector, and a point plus a vector is a point).
Any coordinates system lets us identify the vectors with $\bR^2$, thus making them into a
2-dimensional $\bR$-vector space, and one easily shows that the vector
operations can be defined ``geometrically'' -- independent of the coordinates.
(For multiplication by general real numbers we again need the topology).
Denote the 2-dimensional space of vectors by $V$. By $\End(V)$ we will mean
the space of linear self-maps of $V$.
\bigskip

In so far we had little to do with the primitive equivalence relation of two
pairs of points having the same distance. Now we come to it. By the
requirements from coordinates in Axiom \ref{Ax:Co} any two pairs with the
same vector difference have the same distance, thus we get an equivalence
relation between vectors: {\EM having the same length},
and moreover $v$ and $-v$ always have the same length.
\medskip

Define an {\EM isometry} as an invertible linear self-map $\cU\in\End(V)$
mapping each vector into a vector with same length. The set of isometries
is a group. By the above, $-1$ belongs to this group. (Here and in the sequel
we identify a scalar operator with the scalar).
\medskip

The two remaining axioms deal with isometries. They have a markedly algebraic
flavor, which seems justifiable in view of the above.
\medskip

{\Ax \label{Ax:Ist}
{\EM (Axiom of Isotropy)}. The group of isometries is transitive on a set of
all vectors of the same length, and is also transitive on the set of
$1$-dimensional subspaces of $V$. (That is: for any two vectors of the same
length, or two $1$-dimensional subspaces $\exists$ an isometry mapping one
to the other).
}
\medskip

Instead of the first half of Axiom \ref{Ax:Ist}, one could take the group
of isometries as a primitive notion\NOT{(postulating that they are
(homeomorphisms) that preserves the Middle),} and define vectors to have the
same length iff an isometry maps one to the other.
\medskip

{\Ax \label{Ax:Bd}
{\EM (Axiom of Boundedness)}. The group of isometries is bounded (as a subset
of the $4$-dimensional $\bR$-vector space $\End(V)$).
}
\bigskip

Axiom \ref{Ax:Ist} is related to the empirical/intuitive possibility of motions
(rotations etc.), which is often expressed by congruence axioms.
Axiom \ref{Ax:Bd} postulates that the circle is bounded in a coordinate
system, in spite of the latter extending to infinity in the idealization
which is ``geometry''.
\medskip

Now we shall be able to use the following theorem from algebra/analysis
to obtain that there is a positive-definite quadratic form $Q$ on $V$ such
that vectors $u,v\in V$ have the same length iff $Q(u)=Q(v)$
(thus we have the Theorem of Pythagoras). The resort to such theorem
here seems natural from our point of view. Unfortunately, proving it
requires some mathematical sophistication.

{\Thm \label {Th:Q}
For any bounded group $G\subset\GL(V)$, where $V$ is a 2-dimensional
$\bR$-vector space, there exists a $G$-invariant positive-definite quadratic
form $Q$.}
\bigskip

We give three proofs, differing in the tools used.

\begin{Prff}{\EM Proof 1}
$\BAR{G}$ is a compact group, thus admits a normalized Haar measure $\mu$.
Take any positive-definite
quadratic form $Q_0$ and take as $Q$ the average
$$Q(v)=\int_{g\in\BAR{G}} Q_0(gv)\,d\mu(g).$$
\end{Prff}

This proof works for any finite-dimensional $V$ over $\bR$.

\begin{Prff}{\EM Proof 2}
This again works for any finite-dimensional $V$.
\medskip

Let $W$ be the $\bR$-vector space of quadratic forms on $V$, and $W_+$ the
set of the positive-definite ones (this set is an open convex cone). $G$
acts on $W$ in the canonical way: $(gQ)(v):=Q(g^{-1}v),\; Q\in W$, and
leaves $W_+$ invariant.

Choose a norm $\|\|_0$ on $W$, say the maximum of the absolute values of the
matrix entries with respect to a basis of $V$. Replace $\|\|_0$ by the
$G$-invariant norm
$$\|Q\|:=\sup_{g\in G}\|gQ\|_0.$$
We know that there is a fixed integer $N>0$ such that any subset of $W$ with
$\|\|$-diameter $\le d$ can be covered by at most $N$ sets of $\|\|$-diameter
$\le\frac d2$. if $K\subset W$ is bounded non-empty $G$-invariant convex,
say the convex hull of the orbit of some $Q$, and $\diam(K)\le d$,
then we have a finite set $F\subset K$, of at most $N$ elements, such that
$\forall Q\in K\,\exists Q'\in F\;\|Q-Q'\|\le\frac d2$. This holds, in
particular, for any $Q$ of the form $gQ_0,\;g\in G,\;Q_0:=\dfrac{\sum F}{\# F}$.
Thus $gQ_0$ has distance $\le\frac d2$ from some $Q'\in F$ and distance
$\le d$ from the other members of $F$. This implies $\|gQ_0-Q_0\|<\gamma d$
where $\gamma:=\dfrac{2N-1}{2N}<1$. Since $\|\|$ is $G$-invariant, we have
that the orbit of $Q_0$, hence its convex hull, has diameter $\le\gamma d$.

So we know that any bounded non-empty $G$-invariant convex set $\subset W$
with diameter $\le d$ has a non-empty $G$-invariant convex subset of diameter
$\le\gamma d$. Repeating the process we get an infinite sequence of nested
sets which converges to a $G$-invariant $Q_I\in W$. If we ensure that for any
$Q\in K$ $Q(v)\ge\alpha\|v\|^2$ for some fixed norm $\|\|$ on $V$ and some
fixed $\alpha>0$, then $Q_I$ will be positive-definite.
\end{Prff}

\begin{Prff}{\EM Proof 3}
This is a purely algebraic proof, using $\dim(V)=2$.

There is only a 1-dimensional space of antisymmetric forms on $V$; that is,
a choice of such non-zero form, which we make and denote by
$u\wedge v\;\;u,v\in V$, is possible and unique up to a scalar multiple.
(This follows from $b_1\wedge b_1=0,\;b_2\wedge b_2=0,\;
b_2\wedge b_1=-b_1\wedge b_2$ which any such form must satisfy for a basis
$(b_1,b_2)$, while these formulas indeed give a non-zero antisymmetric form.)

The determinant and trace of the matrix of an $\cA\in\End(V)$ are independent
of the basis (since different bases give similar matrices), hence we may
speak of $\det(\cA)$ and $\tr(\cA)$. The characteristic polynomial of $\cA$
is
$$x^2-\tr(\cA)\cdot x+\det(\cA),\leqno{(1)}$$
its real roots are the real eigenvalues of $\cA$, and plugging $\cA$ in it
gives $0$, by the Cayley-Hamilton Theorem.

For a traceless $\cA\in\End(V)$ (i.e.\ with $\tr(\cA)=0$),
we obtain from the Cayley-Hamilton
Theorem that $\cA^2$ is a scalar, equal to $-\det(\cA)$, and of course to
$\frac12 \tr(\cA^2)$. This scalar gives a quadratic form on the
$3$-dimensional $\bR$-vector space $\{\cA\in\End(v)\,|\,\tr(\cA)=0\}$,
where the corresponding symmetric bilinear form is
$$\LA\cA,\cB\RA:=\frac12\tr(\cA\cB)=\frac12\tr(\cB\cA).$$
Checking the orthogonal basis
$\LP\begin{array}{cc}1&0\\0&-1\end{array}\RP,
  \LP\begin{array}{cc}0&1\\1&0\end{array}\RP,
  \LP\begin{array}{cc}0&1\\-1&0\end{array}\RP$
shows that the signature is $(+,+,-)$. Thus there cannot be two orthogonal
elements with non-positive value of the quadratic form.

If $G\subset\GL(V)$ is a bounded group, then the image of $G$ by $\det$ is a
bounded subgroup of $\bR^\times(=\bR\setminus\{0\})$, therefore
$\det(\cA)=\pm1$ for $\cA\in G$. Also, if $\cA\in G$ has a real eigenvalue
$\lambda$, then $\cA^n\in G$ for integer $n$ and has the eigenvalue
$\lambda^n$, and these must be bounded, therefore $\lambda=\pm1$.

Hence if $\cA\in G$ has determinant $1$, $(1)$ cannot have
a real root different from $\pm1$, which implies $|\tr(\cA)|\le2$. For such
$\cA$, $(1)$ can be written as (recall that $\det(\cA)=1$):
$$\LP x-\frac12\tr(\cA)\RP^2=-\LP1-\LP\frac12\tr(\cA)\RP^2\RP\le0.\leqno{(2)}$$
By Cayley-Hamilton, such $\cA$ can be written as $\frac12\tr(\cA)+\cJ$
where $\cJ$ is traceless with $\cJ^2=-\det(\cJ)$ non-positive, equal to
$-\LP1-\LP\frac12\tr(\cA)\RP^2\RP$.
Also, the product of $\frac12\tr(\cA)\pm\cJ$ is $1$, thus
$\cA^{-1}=\frac12\tr(\cA)-\cJ$.

If we had $\cJ^2=0$ without $\cJ=0$, then $\frac12\tr(\cA)=\pm1$, thus
$\pm\cA=1+\cJ_1$, $\cJ_1^2=0$, $\cJ_1\ne0$ and $(\pm\cA)^n=1+n\cJ_1$
contradicting the boundedness of $G$. Hence if $\cA$ is not $\pm1$ then
$\cJ^2<0$ and $|\tr(\cA)|<2$.

We claim that if $\cA_1=\tau+\cJ_1$ and $\cA_2=\tau+\cJ_2$, $\tau$ scalar,
are elements of $G$ with determinant $1$ and the same trace, $\cJ_1$ and
$\cJ_2$ being traceless with the non-positive square $-(1-\tau^2)$,
then $\cJ_2=\pm\cJ_1$, that is $\cA_2=\cA_1^{\pm1}$.
Indeed, in the above quadratic form on the space of traceless elements of
$\End(V)$, given by the scalar square, $\cJ_1$ and $\cJ_2$ have the same
non-positive square and thus are the sum and difference of the orthogonal
$\frac12(\cJ_1\pm\cJ_2)$. These cannot both have negative square, the
signature being $(+,+,-)$, and none can have square $0$ if $\cJ_1\ne\pm\cJ_2$.
Hence in the latter case they have squares of strictly different signs which
implies
$$|\LA\cJ_1,\cJ_2\RA|=|\frac12\tr(\cJ_1\cJ_2)|>|\LA\cJ_i,\cJ_i\RA|=1-\tau^2.$$
Returning to $\cA_1$ and $\cA_2$ this gives either $\cA_1\cA_2$ or
$\cA_1\cA_2^{-1}$ is an element of $G$ which has half-trace greater than $1$
which we saw above is impossible.

Suppose that we have picked an $\cA\in G$ with determinant $1$ and $\cA$ not
the scalar $\pm1$. We have
$\cA=\frac12\tr(\cA)+\cJ$, $\cA^{-1}=\frac12\tr(\cA)-\cJ$, $\cJ$ is traceless
and $\cJ^2<0$. Consider the symmetric bilinear form on $V$
$$B(u,v):=\frac12\LQ(\cA u)\wedge v+(\cA v)\wedge u\RQ=(\cJ u)\wedge v.
\leqno{(3)}$$
If $b_1$ is a non-zero vector and $b_2=\cJ b_1$, then $b_2$ cannot be
$=\lambda b_1$, $\lambda\in\bR$ because that would imply
$\cJ^2 b_1=\cJ b_2=\lambda \cJ b_1=\lambda b_2=\lambda^2 b_1$.
Therefore $(b_1,b_2)$ is a basis, and $B(b_1,b_1)=b_2\wedge b_1\ne0$.
We have $B(b_1,b_2)=0$ and
$B(b_2,b_2)=(\cJ^2 b_1)\wedge(\cJ b_1)=-(\cJ^2)B(b_1,b_1)$.
So we conclude that $B$ or $-B$ must be positive-definite.

We claim that $B$ is invariant under $G$. Indeed, for $\cB\in G$:
\begin{eqnarray*}
B(\cB u,\cB v)&=&\frac12\LQ(\cA\cB u)\wedge(\cB v)+(cA\cB v)\wedge(\cB u)\RQ=
\\&=&\det(\cB)\frac12\LQ(\cB^{-1}\cA\cB u)\wedge v+
(\cB^{-1}\cA\cB v)\wedge u\RQ
\end{eqnarray*}
But $\cB^{-1}\cA\cB$ is an element of $G$ with determinant $1$ and the same
trace as $\cA$. By the above, it is equal to either $\cA$ or $\cA^{-1}$ and
we find that $B(\cB u,\cB v)$ is one of $\pm B$. Since both are
positive-definite or negative-definite, they are equal.

The theorem is hence proved except when all members of $G$ with determinant
$1$ are scalars. If that is the case, then if there are no elements in $G$
with determinant $-1$ we are done. In any case, for any $\cJ\in G$ with
$\cJ^2=-1$ we can, as above, construct a basis $(b_1,b_2),\;b_2=\cJ b_1$
and the matrix of $\cJ$ in this basis is
$\LP\begin{array}{cc}0&1\\-1&0\end{array}\RP$ with determinant $1$. Hence
such a $\cJ$ is excluded in our case, and all $\cA\in G$ with $\det(\cA)=-1$
must satisfy $\cA^2=1$. Then for any $v\in V$
$v=\frac12(v+\cA v)+\frac12(v-\cA v)$ is a sum of eigenvectors of $\cA$
with eigenvalues $1$ and $-1$, respectively. Since $\det(\cA)=-1$ we
have a $1$-dimensional space of each, i.e.\ $\cA$ has matrix
$\LP\begin{array}{cc}1&0\\0&-1\end{array}\RP$ in some basis. Any other member
of $G$ with determinant $-1$ is a multiple of this $\cA$ by a member of $G$
with determinant $1$, hence the only possible members of $G$ are
$\pm1$, $\pm\cA$ and one easily finds many positive-definite quadratic forms
invariant with respect to these.
\end{Prff}
\medskip

One may remark, following Bourbaki, that after one has the quadratic form on
$V$ that determines equality of length, and thus a corresponding symmetric
bilinear form $\LA\,\RA$, one may define a
$\cJ\in\End V$ by $\LA u,v\RA=(\cJ u)\wedge v$
(see the third proof above) and prove that $\cJ^2$ is a negative scalar,
hence by normalizing the $\wedge$ one may have $\cJ^2=-1$ (which determines
$\cJ$ up to sign). This turns $V$ in a canonical way into a $1$-dimensional
complex vector space (didactically, the complex numbers may be defined by
our ``geometric'' $\cJ$ -- in this approach every plane (say, in $3$-space)
has, strictly speaking, its own ``complex numbers''.) Using this complex
structure to do plane geometry is very fruitful. For example, angles (with
their trigonometry) can be easily treated in a rigorous way.

\section{Space Geometry}

To axiomatize the Euclidean geometry of $n$-space, one may start from a set
{\EM Space} of points with exactly analogous primitive notions, replace
$\bR^2$ by $\bR^n$ in the Axiom of Coordinates, and postulate the Axioms of
Isotropy and of Boundedness for every sub-$2$-plane of Space
(or alternatively for Space itself). In passing from the existence of an
equality-of-length -- determining quadratic form on every plane to the
existence of one such form for Space, one may use the well-known

{\Thm
Let $V$ be an (not necessarily finite-dimensional) $\bR$-vector space,
$\dim(V)\ge2$, and let $U\subset V$. If for every $2$-dimensional subspace
(=plane) $V'\subset V$ there exists a positive-definite quadratic form $Q'$ on $V'$
such that $U\cap V'=\{v\in V'\,|\,Q'(v)=1\}$ then there exists a
positive-definite quadratic form $Q$ on the whole $V$ such that
$U=\{v\in V\,|\,Q(v)=1\}$.
}
\begin{Prf}
It is clear that the $Q'$ are unique for each plane, and that they agree
on intersections. Hence they define a function $Q$ on $V$, positive on
$V\setminus{0}$ and it remains to prove that $Q$ is quadratic, i.e.\ comes
from a symmetric bilinear form on $V$. That form must be:
$$\LA u,v\RA=\frac12(Q(u+v)-Q(u)-Q(v))\quad u,v\in V$$
and we have to prove that this is bilinear. Since we have
$\LA\lambda u,v\RA=\lambda\LA u,v\RA$, $Q$ being quadratic on the plane
containing $u$ and $v$, it remains to prove biadditivity.
As we clearly have $\LA u,0\RA=0$,
biadditivity will follow if we prove that $\LA u,v\RA+\LA u,w\RA$ depends
only on $u$ and $v+w$. This obtains from the following calculation
(where one uses the parallelogram
equality $2Q(a)+2Q(b)=Q(a+b)+Q(a-b),\;\;a,b\in V$, which holds since $Q$
is quadratic in the plane containing $a$ and $b$):
\begin{eqnarray*}
&&2\LP\LA u,v\RA+\LA u,w\RA\RP=Q(u+v)-Q(u)-Q(v)+Q(u+w)-Q(u)-Q(w)=\\
&&=\frac12\LP Q(2u+v+w)+Q(v-w)\RP-2Q(u)-\frac12\LP Q(v+w)+Q(v-w)\RP=\\
&&=\frac12\LP Q(2u+v+w)-Q(v+w)\RP-2Q(u)
\end{eqnarray*}
\end{Prf}

\end{document}